\documentclass[11pt]{article}
\usepackage{amsfonts, amssymb, amsmath}
\usepackage{enumerate}
\usepackage{esint}
\newtheorem{thm}{Theorem}[section]

\newtheorem{cor}[thm]{Corollary}
\newtheorem{lem}[thm]{Lemma}
\newtheorem{prop}[thm]{Proposition}
\newtheorem{defn}[thm]{Definition}
\newtheorem{rem}[thm]{Remark}

\newcommand{\f}{\frac}

\newcommand{\vc}{\infty}

\newcommand{\RR}{\mathbb{R}}

\newcommand{\pa}{\partial}
\begin{document}
\title{Weighted norm inequalities for pseudo-differential operators and their commutators }
\author{The Anh Bui\thanks{Department of Mathematics, Macquarie University, NSW 2109, Australia and Department of Mathematics, University of Pedagogy, HoChiMinh City, Vietnam. \newline{Email: the.bui@mq.edu.au and bt\_anh80@yahoo.com}  \newline
{\it {\rm 2000} Mathematics Subject Classification:} 42B20, 35S05, 47G30.
\newline
{\it Key words:} Weighted norm inequality; Pseudo-differential operator; Commutator.}}

\date{}

\maketitle

\begin{abstract}
This paper is dedicated to study weighted $L^p$ inequalities for pseudo-differential operators with amplitudes and their
commutators by using the new class of weights $A_p^\vc$ and the new BMO function space BMO$_\vc$ which are larger than the Muckenhoupt class of weights $A_p$ and classical BMO space BMO, respectively. The obtained results therefore improve substantially some well-known results.
\end{abstract}

 \tableofcontents
\section{Introduction and the main results}

For $f \in C_0^\vc(\RR^n)$ a pseudo-differential operator given formally by
$$
T_af(x)=\f{1}{(2\pi)^n}\int_{\RR^n}\int_{\RR^n}a(x,y,\xi)e^{i\langle x-y,\xi \rangle}f(y)dyd\xi,
$$
where the amplitude $a$ satisfies certain growth conditions. The boundedness of pseudo-differential operators has been studied extensively by many mathematicians, see for example \cite{AH, CT, H1, H2, MRS, N, Y} and the references therein. One of the most interesting problems is studying the weighted norm inequalities for pseudo-differential operators and their commutators with BMO function, see for example \cite{M, MRS, N,T,Y}.

\medskip

In this paper we consider the following classes of symbols and amplitudes $a$ (in what follows we set $\langle x \rangle=(1+|x|^2)^{1/2}$):

\begin{defn}\label{defn1}
Let $a:\RR^n\times \RR^n\times \RR^n \rightarrow \RR^n$ and $m\in \RR, \rho\in [0,1]$ and $\delta\in [0,1]$.
\begin{enumerate}[(a)]
\item We say $a\in A_{\rho,\delta}^m$ when for each triple of multi-indices $\alpha, \beta$ and $\gamma$ there exists a constant $C$ such that
$$
|\pa_\xi^\alpha\pa_x^\beta\pa_y^\gamma a(x,y,\xi)|\leq C\langle \xi \rangle^{m-\rho|\alpha|+\delta|\beta+\gamma|}.
$$
\item We say $a\in L^\vc A_{\rho,\delta}^m$ when for each triple of multi-indices $\alpha, \beta$ and $\gamma$ there exists a constant $C$ such that
$$
\|\pa_\xi^\alpha\pa_y^\beta a(\cdot,y,\xi)\|_{L^\vc}\leq C\langle \xi \rangle^{m-\rho|\alpha|+\delta|\beta|}.
$$

\end{enumerate}

\end{defn}

\begin{defn}\label{defn2}
Let $a:\RR^n\times \RR^n \rightarrow \RR^n$ and $m\in \RR, \rho\in [0,1]$ and $\delta\in [0,1]$. \\

\begin{enumerate}[(a)]
\item We say $a\in S_{\rho,\delta}^m$ when for each pair of multi-indices $\alpha$ and $\beta$ there exists a constant $C$ such that
$$
|\pa_\xi^\alpha\pa_x^\beta a(x,\xi)|\leq C\langle \xi \rangle^{m-\rho|\alpha|+\delta|\beta|}.
$$

\item We say $a\in L^\vc S_{\rho}^m$ when for each multi-indices $\alpha$ there exists a constant $C$ such that
$$
\|\pa_\xi^\alpha  a(\cdot,\xi)\|_{L^\vc}\leq C\langle \xi \rangle^{m-\rho|\alpha|}.
$$
\end{enumerate}
\end{defn}
It is easy to see that $S_{\rho,\delta}^m\subset A_{\rho,\delta}^m$, $L^\vc S^m_\rho \subset L^\vc A_{\rho,\delta}^m$, $S_{\rho,\delta}^m\subset L^\vc S_{\rho}^m$ and $A_{\rho,\delta}^m\subset L^\vc A_{\rho,\delta}^m$. The classes $A_{\rho,\delta}^m$ and $S_{\rho,\delta}^m$ were studied in \cite{M, H1}. For further information about these two classes, we refer the reader to for example \cite{H1, St}. The class $L^\vc S_{\rho}^m$ were introduced by \cite{MRS} and it is the natural generalization of the class $S_{\rho,\delta}^m$. This class is much rougher than that considered in \cite{N, Y}.

 The aim of this paper is to study the weighted norm inequalities for pseudo-differential operators $T_a$ and their commutators by using the new BMO functions and the new class of weights. Firstly, we would like to give brief definitions on the new class of weights and the new BMO function space (we refer to Section 2 for details):

 \medskip

 The new classes of weights $A^\vc_{p}=\cup_{\theta>0}A^{\theta}_{p}$ for $p\geq 1$, where $A^{\theta}_{p}$, $\theta\geq 0$, is the set of those weights satisfying
\begin{equation}\label{classofnewweights}
\Big(\int_Bw\Big)^{1/p}\Big(\int_Bw^{-\f{1}{p-1}}\Big)^{1/p'}\leq C|B|(1+r_B)^{\theta}
\end{equation}
for all ball $B=B(x_B,r_B)$. We denote $A^\vc_{\vc}=\cup_{p\geq 1}A^{\vc}_{p}$. It is easy to see that the new class $A_p^\vc$ is larger than the Muckenhoupt class $A_p$.

\medskip

The new BMO space $BMO_\theta$ with $\theta\geq 0$ is defined as a set of all locally integrable functions $b$ satisfying
\begin{equation}\label{eq1-intro}
\f{1}{|B|}\int_B|b(y)-b_B|dy\leq C(1+r_B)^\theta
\end{equation}
where $B=B(x_B,r_B)$ and $b_B=\f{1}{|B|}\int_B b$. A norm for $b\in BMO_\theta$, denoted by $\|b\|_\theta$, is given by the infimum of the constants satisfying (\ref{eq1-intro}). Clearly $BMO_{\theta_1}\subset BMO_{\theta_2}$ for $\theta_1\leq \theta_2$ and $BMO_0=BMO$. We define $BMO_\vc=\cup_{\theta>0}BMO_\theta$.

\medskip

Our main result is the following theorem.

\begin{thm}\label{mainresult}
Let $a\in L^\vc A^m_{\rho,\delta}$ with $m<n(\rho-1)$ or $a\in L^\vc A^0_{1,\delta}, \delta\in [0,1]$. If $T_a$ is bounded on $L^p$ for all $1<p<\vc$, then

\begin{enumerate}[(a)]
\item $T_a$ is bounded on $L^p(w)$ for $1<p<\vc$ and $w\in A_p^\vc$;

\item For any $b\in BMO_\vc$, the commutator $[b, T_a]$ bounded on $L^p(w)$ for $1<p<\vc$ and $w\in A_p^\vc$.
\end{enumerate}
\end{thm}

We would like to specify some applications of Theorem \ref{mainresult}:

\medskip

In \cite{M}, the author study the weighted $L^p$ inequalities of $T_a$ when the symbol $a$ belongs to the class $S^0_{1, \delta}\subset L^\vc A^0_{1,\delta}$ with $\delta\in (0,1)$. It was proved that $T_a$ is bounded on $L^p(w)$ for $1<p<\vc$, $w\in A_p$. Recently, the author in \cite{T} showed that $T_a$ and its commutator with BMO function $[b, T_a]$ is bounded on $L^p(w)$   for $1<p<\vc$ and $w\in A^\vc_p$ by the different approach. Here, by using Theorem \ref{mainresult}, we not only re-obtain the boundedness of $T_a$ on $L^p(w)$ for $1<p<\vc$ and $w\in A^\vc_p$, but also obtain the new result on the boundedness of its commutator with BMO$_\vc$ functions.

\begin{cor}\label{cor1}
Let $a\in S^0_{1, \delta}\subset L^\vc A^0_{1,\delta}, 0<\delta<1$. Then we have

(i) $T_a$ is bounded on $L^p(w)$   for $1<p<\vc$ and $w\in A^\vc_p$;

(ii) For each $b\in BMO_\vc$, the commutator $[b, T_a]$ is bounded on $L^p(w)$   for $1<p<\vc$ and $w\in A^\vc_p$.
\end{cor}

Now we consider the class $L^\vc S^m_\rho$. If $a\in L^\vc S^m_\rho$ with $\rho\in [0,1]$ and $m<n(\rho-1)$, then the authors in \cite{MRS} proved that the pseudo-differential operator $T_a$ and its commutators with BMO functions $[b, T_a]$ are bounded on $L^p(w)$ for $1<p<\vc$ and $w\in A_p$, see Theorem 3.3 and 4.5 in \cite{MRS}. So, Theorem \ref{mainresult} leads us to the following result.

\begin{cor}\label{cor2}
Let $a\in L^\vc S^m_\rho$ with $\rho\in [0,1]$ and $m<n(\rho-1)$. Then we have

(i) $T_a$ is bounded on $L^p(w)$   for $1<p<\vc$ and $w\in A^\vc_p$;

(ii) For each $b\in BMO_\vc$, the commutator $[b, T_a]$ is bounded on $L^p(w)$   for $1<p<\vc$ and $w\in A^\vc_p$.
\end{cor}

 It was proved in \cite[Theorem 3.7]{MRS} that if $a\in L^\vc A^m_{\rho,\delta}$ with $0\leq \rho\leq 1$ and $m<n(\rho-1)$,  then $T_a$ and $[b, T_a]$ are bounded on $L^p(w)$ for $1<p<\vc$ and $w\in A_p$ with $b\in BMO$. Therefore, in the light of Theorem \ref{mainresult}, we have:
\begin{cor}\label{cor4}
Let $a\in L^\vc A^m_{\rho,\delta}$ with $0\leq \rho\leq 1$ and $m<n(\rho-1)$. Then we have

(i) $T_a$ is bounded on $L^p(w)$   for $1<p<\vc$ and $w\in A^\vc_p$;

(ii) For each $b\in BMO_\vc$, the commutator $[b, T_a]$ is bounded on $L^p(w)$   for $1<p<\vc$ and $w\in A^\vc_p$.
\end{cor}

For the smooth amplitude, we have the following result.
\begin{cor}\label{cor6}
Let $a\in A^{n(\rho-1)}_{\rho,\delta}$ with $0<\rho\leq 1, 0\leq \delta<1$.  Then we have

(i) $T_a$ is bounded on $L^p(w)$   for $1<p<\vc$ and $w\in A^\vc_p$;

(ii) For each $b\in BMO_\vc$, the commutator $[b, T_a]$ is bounded on $L^p(w)$   for $1<p<\vc$ and $w\in A^\vc_p$.
\end{cor}
\emph{Proof:} The remark in \cite[p. 11]{AH} tells us that $T_a$ is bounded on $L^p$ for $1<p<\vc$. Thanks to Theorem \ref{mainresult}, we conclude that $T_a$ and $[b, T_a], b\in BMO_\vc$ are bounded on $L^p(w)$ for $1<p<\vc$ and $w\in A^\vc_p$.

\medskip


The outline of the paper is as follows. In Section 2, we first recall some definitions of the new class of weights $A_p^\vc$ and the new BMO function spaces $BMO_\vc$. Then we also review some basic properties concerning on $A_p^\vc$ and $BMO_\vc$. Section 3 represents some kernel estimates for the pseudo-differential operator $T_a$. The proof of the main result will be given in Section 4.
\section{Preliminaries}

To simplify notation, we will often just use $B$ for $B(x_B, r_B)$ and $|E|$ for the measure of $E$ for any measurable subset $E\subset \RR^n$.
Also given $\lambda > 0$, we will write $\lambda B$ for the
$\lambda$-dilated ball, which is the ball with the same center as
$B$ and with radius $r_{\lambda B} = \lambda r_B$. For each ball
$B\subset \RR^n$ we set
$$
S_0(B)=B \ \text{and} \ S_j(B) = 2^jB\backslash 2^{j-1}B \
\text{for} \ j\in \mathbb{N}.
$$

\subsection{The new class of weights and new BMO function spaces}
Recently, in \cite{BHS2}, a new class of weights  associated to Schr\"odinger operators $L:=-\Delta +V$ where the potential $V \in RH_{n/2}$, the reverse H\"older class has been introduced. According to \cite{BHS2}, the authors defined the new classes of weights $A^L_{p}=\cup_{\theta\geq 0}A^{L,\theta}_{p} $ for $p\geq 1$ , where $A^{L,\theta}_{p}$, $\theta\geq 0$, is the set of those weights satisfying
\begin{equation}\label{classofnewweights}
\Big(\int_Bw\Big)^{1/p}\Big(\int_Bw^{-\f{1}{p-1}}\Big)^{1/p'}\leq C|B|\Big(1+\f{r}{\rho(x)}\Big)^{\theta}
\end{equation}
for all ball $B=B(x,r)$. We denote $A^L_{\vc}=\cup_{p\geq 1}A^{L}_{p}$ where the critical radius function $\rho(\cdot)$ is defined by
\begin{equation}
\rho(x)=\sup\Big\{r>0:\f{1}{r^{n-2}}\int_{B(x,r)}V\leq 1\Big\}, \ \ x\in \mathbb{R}^n.
\end{equation}

In this paper, we consider the particular case when $\rho(\cdot)\equiv 1$. In this situation the new classes of weights is defined by $A^\vc_{p}=\cup_{\theta\geq 0}A^{\theta}_{p}$ for $p\geq 1$, where $A^{\theta}_{p}, \theta\geq 0,$ is the set of those weights satisfying
\begin{equation}\label{classofnewweights}
\Big(\int_Bw\Big)^{1/p}\Big(\int_Bw^{-\f{1}{p-1}}\Big)^{1/p'}\leq C|B|(1+r_B)^{\theta}
\end{equation}
for all ball $B=B(x_B,r_B)$. We denote $A^\vc_{\vc}=\cup_{p\geq 1}A^{\vc}_{p}$.

It is easy to see that the new class $A_p^\vc$ is larger than the Muckenhoupt class $A_p$.
The following properties hold for the new classes $A^\vc_p$, see \cite[Proposition 5]{BHS2}.
\begin{prop}\label{property-Avc}
 The following statements hold:

i) $A_{p}^\vc\subset A_q^\vc$ for $1\leq p\leq q<\vc$.

ii) If $w\in A^\vc_p$ with $p> 1$ then there exists $\epsilon >0$ such that $w\in A_{p-\epsilon}^\vc$. Consequently, $A_p^\vc=\cup_{q<p}A_q^\vc$.
\end{prop}

Similarly, by adapting the ideas to \cite{BHS1}, the new BMO space $BMO_\theta$ with $\theta\geq 0$ is defined as a set of all locally integrable functions $b$ satisfying
\begin{equation}\label{eq1-intro}
\f{1}{|B|}\int_B|b(y)-b_B|dy\leq C(1+r_B)^\theta
\end{equation}
where $B=B(x_B,r_B)$ and $b_B=\f{1}{|B|}\int_B b$. A norm for $b\in BMO_\theta$, denoted by $\|b\|_\theta$, is given by the infimum of the constants satisfying (\ref{eq1-intro}). Clearly $BMO_{\theta_1}\subset BMO_{\theta_2}$ for $\theta_1\leq \theta_2$ and $BMO_0=BMO$. We define $BMO_\vc=\cup_{\theta>0}BMO_\theta$.

The following result can be considered to be a variant of John-Nirenberg inequality for the spaces $BMO_L^\theta$.
\begin{prop}\label{JNforBMOL}
Let $\theta>0, s\geq 1$. If $b\in BMO_L^\theta$ then for all $B=(x_0, r)$

i) $$
\Big(\f{1}{|B|}\int_B|b(y)-b_B|^sdx\Big)^{1/s}\lesssim \|b\|_{\theta}(1+r_B)^{\theta};
$$

ii) $$
\Big(\f{1}{|2^kB|}\int_{2^kB}|b(y)-b_B|dx\Big)^{1/s}\lesssim \|b\|_{\theta} k(1+2^kr_B)^{\theta}
$$
for all $k\in \mathbb{N}$.

\end{prop}

The proof is similar (even easier) to Lemma 1 and Proposition 3 in \cite{BHS2} and hence we omit details.

\subsection{Weighted estimates for some localized operators}

A ball of the form $B(x_B,r_B)$ is called {\it a critical ball} if $r_B =1$. We have the following result.

\begin{prop}\label{coveringlemma}
There exists a sequence of points $x_j, j\geq 1$ in $\mathbb{R}^n$ so that the family of critical balls $\{Q_j\}_j$ where $Q_j:=B(x_j, 1)$, $j\geq 1$ satisfies

(i) $\cup_j Q_j = \mathbb{R}^n$.

(ii) There exists a constant $C$ such that for any $\sigma>1$, $\sum_j \chi_{\sigma Q_j}\leq C\sigma^{n}$.
\end{prop}

Note that the more general version of Proposition \ref{coveringlemma} is obtained by \cite{DZ}. However, in our particular situation, for convenience, we would like to give a simple proof of this proposition.

\medskip

\emph{Proof:} Let us consider the family of balls $\{B(x,\f{1}{5}): x\in \RR^n\}$. Using Vitali covering lemma, we can pick the subfamily of balls $\{B_j:=B(x_j,\f{1}{5}): j\geq 1\}$ so that $\{Q_j\}_j$ is pairwise disjoint and $\RR^n\subset \cup_j Q_j$ where $Q_j=5B_j= B(x_j, 1)$. This gives (i).

To prove (ii), pick any $x\in \RR^n$. Let $\mathfrak{I}$ be the set of all indices $j$ so that $x\in \sigma Q_j$. Note that if $x\in \sigma Q_j$ then $\sigma Q_j\subset B(x, 2\sigma)$. Therefore, $B(x_j, \f{1}{5})\subset B(x, 2\sigma)$ for all $j\in \mathfrak{I}$. Since $\{B(x_j, \f{1}{5})\}_{j\in \mathfrak{I}}$ is pairwise disjoint, $\sum_{j\in \mathfrak{I}}|B(x_j, \f{1}{5})|\leq |B(x, 2\sigma)|$. This is equivalent to that $|\mathfrak{I}|/5^n\leq C\sigma^n$. Hence, $|\mathfrak{I}|\leq C\sigma^n$. This completes our proofs.

\medskip

We consider the following maximal functions for $g\in L^1_{{\rm loc}}(\mathbb{R}^n)$ and $x\in \mathbb{R}^n$
$$
M_{{\rm loc}, \alpha}g(x)=\sup_{x\in B\in \mathcal{B}_{\alpha}}\f{1}{|B|}\int_B|g|,
$$
$$
M^\sharp_{{\rm loc}, \alpha}g(x)=\sup_{x\in B\in \mathcal{B}_{\alpha}}\f{1}{|B|}\int_B|g-g_B|,
$$
where $\mathcal{B}_{\alpha}=\{B(y,r): y\in \mathbb{R}^n \ \text{and} \ r\leq \alpha \}$.

Also, given a ball $Q$, we define the following maximal functions for $g\in L^1_{{\rm loc}}(\mathbb{R}^n)$ and $x\in Q$
$$
M_{Q}g(x)=\sup_{x\in B\in \mathcal{F}(Q)}\f{1}{|B\cap Q|}\int_{B\cap Q}|g|,
$$
$$
M^\sharp_{Q}g(x)=\sup_{x\in B\in F(Q)}\f{1}{|B\cap Q|}\int_{B\cap Q}|g-g_{B\cap Q}|,
$$
where $\mathcal{F}(Q)=\{B(y,r): y\in Q, r>0 \}$.

We have the following lemma.
\begin{lem}\label{FSinequalityversion}
For $1<p<\vc$, then there exists $\beta$ such that if $\{Q_k\}_{k}$ is a sequence of balls as in Proposition \ref{coveringlemma} then
$$
\int_{\mathbb{R}^n}|M_{{\rm loc},\beta}g(x)|^pw(x)dx \lesssim \int_{\mathbb{R}^n}|M^\sharp_{{\rm loc},4}g(x)|^pw(x)dx+\sum_{k}w(Q_k)\Big(\f{1}{|2Q_k|}\int_{2Q_k}|g|\Big)^p
$$
for all $g\in L^1_{{\rm loc}}(\mathbb{R}^n)$ and $w\in A^\vc_{\vc}.$
\end{lem}
\emph{Proof:} We adapt the argument in \cite[Lemma 2]{BHS1} to our present situation.

Taking $\beta=1/2$, by Lemma \ref{coveringlemma}, we have
$$
\int_{\mathbb{R}^n}|M_{{\rm loc},\f{1}{2}}g(x)|^pw(x)dx \leq C\sum_{k}\int_{Q_k}|M_{{\rm loc},\f{1}{2}}g(x)|^pw(x)dx.
$$
It can be verified that for $x\in Q_k$, $M_{{\rm loc},\f{1}{2}}g(x)\leq M_{2Q_k}(g\chi_{2Q_k})$. Note that since $g\chi_{2Q_k}$ is supported in $2Q_k$, operators $M_{2Q_k}$ and $M_{2Q_k}^\sharp$ are
Hardy-Littlewood and sharp maximal functions defined in $2Q_k$
viewed as a space of homogeneous type with the Euclidean metric and the Lebesgues
measure restricted to $2Q_k$. Moreover, by definition of $A_\vc^\vc$, if $w\in A_\vc^\vc$ then $w\in A_\vc(2Q_k)$, where $A_\vc(2Q_k)=\cup_{p\geq 1}A_p(2Q_k)$ and $A_p(2Q_k)$ is the class of Muckenhoupt weights on the spaces of homogeneous type $2Q_k$. Moreover, due to \cite[Lemma 5]{BHS2}, $[w]_{A_\vc(2Q_k)}\leq C$ for all $k\geq 1$. Therefore, using Proposition 3.4 in \cite{PS} gives
\begin{equation*}
\begin{aligned}
\int_{\mathbb{R}^n}&|M_{{\rm loc},\f{1}{2}}g(x)|^pw(x)dx\\
 &\leq C\sum_{k}\int_{Q_k}|M_{{\rm loc},\f{1}{2}}g(x)|^pw(x)dx\\
& \leq C\sum_{k}\int_{Q_k}|M_{2Q_k}(g\chi_{2Q_k})(x)|^pw(x)dx\\
& \leq C\sum_{k}\int_{2Q_k}|M^\sharp_{2Q_k}(g\chi_{2Q_k})(x)|^pw(x)dx+C\sum_k w(2Q_k)\Big(\f{1}{|2Q_k|}\int_{2Q_k}|g(x)|w(x)dx\Big)^p.
\end{aligned}
\end{equation*}
To complete the proof, we need only to check that $M^\sharp_{2Q_k}(g\chi_{2Q_k})(x)\leq C M^\sharp_{{\rm loc}, 4}(g)(x), x\in 2Q_k$. We have
$$
M^\sharp_{{\rm loc}, 4}(g)(x)=\sup_{B\in F(2Q_k): B\ni x}\f{1}{|B\cap 2Q_k|}\int_{B\cap 2Q_k}|f-f_{B\cap 2Q_k}|.
$$
If $r_B\geq 4$, due to $r_{2Q_k}=2$, $2Q_k\subset B$. Hence, in this situation, we have
$$
\f{1}{|B\cap 2Q_k|}\int_{B\cap 2Q_k}|f-f_{B\cap 2Q_k}|=\f{1}{|Q_k|}\int_{2Q_k}|f-f_{2Q_k}|\leq M^\sharp_{{\rm loc}, 4}(g)(x).
$$
Otherwise, if $r_B<4$, it is obvious that $|B\cap 2Q_k|\approx |B|$. So we have
\begin{equation*}
\begin{aligned}
\f{1}{|B\cap 2Q_k|}\int_{B\cap 2Q_k}|f-f_{B\cap 2Q_k}|&\leq 2\f{1}{|B\cap 2Q_k|}\int_{B\cap 2Q_k}|f-f_{B}|\\
&\leq C\f{1}{|B|}\int_{B}|f-f_{B}|\leq CM^\sharp_{{\rm loc}, 4}(g)(x).
\end{aligned}
\end{equation*}
This completes our proof.\\

\medskip

Throughout this paper, we always assume that $N$ is a sufficiently large number and different from line to line. For $\kappa\geq 1$ and $p\geq 1$, we define the following functions for $g\in L^1_{{\rm loc}}(\mathbb{R}^n)$  and $x\in \mathbb{R}^n$
$$
G_{\kappa, p}f(x)=\sup_{Q\ni x; Q \ {\rm is \ critical}}\sum_{k=0}^\vc 2^{-Nk}\Big(\f{1}{|2^k\widehat{Q}|}\int_{2^k\widehat{Q}}|f(z)|^pdz\Big)^{1/p}
$$
where $\widehat{Q}=\kappa Q$.

When $\kappa=1$, we write $G_p$ instead of $G_{1,p}$. The following result gives the  weighted estimates for $G_{\kappa,p}$.
\begin{prop}\label{weighted estiamtes for G}
Let $p>s>1$ and $w\in A_{p/s}^\theta$, $\theta\geq 0$. Then we have
$$
\|G_{\kappa, s}f\|_{L^p(w)}\lesssim \|f\|_{L^p(w)}.
$$
\end{prop}
Without the loss of generality, we assume that $\kappa=1$. Assume that $Q=B(x_0,1)$. For $x\in Q$, $Q\subset 2B(x,1)$. This implies that
$$
G_pf(x)\leq  C \sum_{k=0}^\vc 2^{-N k}\Big(\f{1}{|2^kB(x,1)|}\int_{B_k(x,1))}|f(z)|^sdz\Big)^{1/s}
$$
where $B_k(x,1)=B(x,2^{k+1})$.

Let $\{Q_j\}$ be the family of critical balls given by Proposition \ref{coveringlemma}. Note that if $x\in Q_j$, $B_k(x,1)\subset Q^k_j$ where $Q^k_j= 2^{k+2}Q_j$. These estimates and H\"older inequalities give
\begin{equation}\label{eq1-G}
\begin{aligned}
\|G_pf\|_{L^p(w)}&\leq  C  \sum_{k=0}^\vc 2^{-N k}\Big(\sum_j\int_{Q_j}\Big(\f{1}{|2^kB(x,1)|}\int_{B_k(x,1)}|f(z)|^sdz\Big)^{p/s}w(x)dx\Big)^{1/p}\\
&\leq  C  \sum_{k=0}^\vc 2^{-N k}\Big(\sum_j\int_{Q_j}\Big(\f{1}{|2^kQ_j|}\int_{Q_j^k}|f(z)|^sdz\Big)^{p/s}w(x)dx\Big)^{1/p}\\
&\leq  C  \sum_{k=0}^\vc 2^{-N k}\Big(\sum_j\f{w(Q_j)}{|2^kQ_j|^{p/s}}\Big(\int_{Q_j^k}|f(z)|^sdz\Big)^{p/s}\Big)^{1/p}\\
&\leq  C  \sum_{k=0}^\vc 2^{-N k}\Big(\sum_j \f{w(Q^k_j)}{|2^kQ_j|^{p/s}} \Big(\int_{Q_j^k}w^{-\f{(p/s)'}{p/s}}\Big)^{\f{p/s}{(p/s)'}}\Big(\int_{Q_j^k}|f(z)|^pw(z)dz\Big)\Big)^{1/p}.
\end{aligned}
\end{equation}
Since $w\in A_{p/s}^{\theta}$, by definition of the classes $A_p^{\theta}$, we have
$$
w(Q^k_j) \Big(\int_{Q_j^k}w^{-\f{(p/s)'}{p/s}}\Big)^{\f{p/s}{(p/s)'}}\leq C|Q_j^k|^{p/s}2^{k\theta\times (p/s)}.
$$

This together with (\ref{eq1-G}) gives
\begin{equation*}
\begin{aligned}
\|G_sf\|_{L^p(w)}&\leq  C \sum_{k}2^{-k(N-\theta/s)}\Big(\sum_j\int_{Q_j^k}|f(z)|^pw(z)dz\Big)^{1/p}\\
&\leq  C \sum_{k}2^{-k(N-\theta/s-n/p)}\|f\|_{L^p(w)}\\
&\leq  C\|f\|_{L^p(w)}.
\end{aligned}
\end{equation*}
This completes our proof.

\medskip

For a family of balls $\{Q_k\}_k$ given by Proposition \ref{coveringlemma}, we define the operator $\widetilde{M}_{s}, s\geq 1,$ as follows
\begin{equation}\label{defnoftidleM}
\widetilde{M}_sf=\sum_k\chi_{Q_k}M_s(f\chi_{\widetilde{Q}_k})
\end{equation}
where $\widetilde{Q}_j= 8 Q_j$ and $M_sf=M(|f|^s)^{1/s}$ with $M$ is the Hardy-Littlewood maximal function. We have the following result.
\begin{prop}\label{rem1}
If $p>s>1$ and $w\in A^\theta_{p/s},\theta>0$, then $\widetilde{M}_s$ is bounded on $L^p(w)$.
\end{prop}
\emph{Proof:}
We have
$$
\int_{\RR^n}|\widetilde{M}_sf(x)|^pw(x)dx=\sum_{j}\int_{Q_j}|M_s(f\chi_{\widetilde{Q}_k})|^pw(x)dx.
$$
For each $k$, if we consider $\widetilde{Q}_k$ as a space of homogeneous type with the Euclidean metric and the
Lebesgues measure restricted to $\widetilde{Q}_k$, then $w\in A_{p/s}(\widetilde{Q}_k)$. Moreover, it can be verified that
$$
\|M_s(f\chi_{\widetilde{Q}_k})\|_{L^p(w, \widetilde{Q}_k)}\leq C\|f\|_{L^p(w, \widetilde{Q}_k)}
$$
and the constant $C$ is independent of $k$.

Therefore, by (ii) of Lemma \ref{coveringlemma},
$$
\begin{aligned}
\int_{\RR^n}|\widetilde{M}_sf(x)|^pw(x)dx&\leq C\sum_{j}\int_{\widetilde{Q}_k}|f(x)|^pw(x)dx\\
&\leq C\|f\|_{L^p(w)}^p.
\end{aligned}
$$
This completes our proof.

\section{Some kernel estimates}

Let $\varphi_0:\RR^n \rightarrow \RR$ be a smooth radial function which is equal to $1$ on the unit ball centered at origin and supported on its concentric double. Set $\varphi(\xi)=\varphi_0(\xi)-\varphi_0(2\xi)$ and $\varphi_k(\xi)=\varphi(2^{-k}\xi)$. Then, we have
$$
\sum_{k=0}^\vc\varphi_k(\xi)=1 \ \text{for all $\xi\in \RR^n$}
$$
and supp $\varphi_k\subset \{\xi: 2^{k-1}\leq |\xi|\leq 2^{k+1}\}$ for all $k\geq 1$. Moreover, for  any multi-index $\alpha$ and $N\geq 0$, we have
$$
|\pa^\alpha_\xi \varphi_k(\xi)|\leq c_\alpha 2^{-k|\alpha|}.
$$

\begin{lem}\label{lem1}
Let $a\in L^\vc A^m_{\rho,\delta}$ with $m\in \RR, \rho\in [0,1]$  and $\delta\in [0,1]$. Let $a_k(x,y,\xi)=a(x,y,\xi)\varphi_k(\xi)$ for $k\geq 0$.
\begin{enumerate}[(a)]
\item For each $\ell \geq 0$,
$$
|z|^\ell \Big|\int a_k(x,y,\xi)e^{i\langle z,\xi\rangle}d\xi\Big|\leq C2^{k(n+m-\rho\ell)}.
$$

\item If $a\in L^\vc A^m_{\rho,\delta}$ with $m<n(\rho-1)$ and $\rho,\delta\in [0,1]$,  then for each $N>0, $there exist $\epsilon , \epsilon'>0$ so that for any ball $B\subset \RR^n$, $y, \overline{y} \in B,$ and $x \in S_j(B), j\geq 2$ so that
$$
\Big|\int a_k(x,y,\xi)e^{i\langle x-y,\xi\rangle}-a_k(x,\overline{y},\xi)e^{i\langle x-\overline{y},\xi\rangle}d\xi\Big|\leq C2^{-j\epsilon}(2^jr_B)^{-n}\min\{1, (2^jr_B)^{-N}\} 2^{-k\epsilon'}.
$$

\item Particularly, if $a\in L^\vc A^{0}_{1,\delta}$, $\delta\in [0,1]$,
then there exist $\epsilon , \epsilon'>0$ so that for any ball $B\subset \RR^n$, $y, \overline{y} \in B,$ and $x \in S_j(B), j\geq 2$ so that
$$
\Big|\int a_k(x,y,\xi)e^{i\langle x-y,\xi\rangle}-a_k(x,\overline{y},\xi)e^{i\langle x-\overline{y},\xi\rangle}d\xi\Big|\leq C2^{-j\epsilon}(2^jr_B)^{-n}\min\{1, (2^jr_B)^{-N}\}(2^{k}r_B)^{\epsilon'}
$$
as long as $2^{k}r_B\leq 1$; and
$$
\Big|\int a_k(x,y,\xi)e^{i\langle x-y,\xi\rangle}-a_k(x,\overline{y},\xi)e^{i\langle x-\overline{y},\xi\rangle}d\xi\Big|\leq C2^{-j\epsilon}(2^jr_B)^{-n}\min\{1, (2^jr_B)^{-N}\}(2^{k}r_B)^{-\epsilon'}
$$
as long as $2^{k}r_B> 1$.
\end{enumerate}
\end{lem}
\emph{Proof:} We refer Lemma 3.1 in \cite{MRS} for the proof of (a).

\medskip

(b)
We first note that since $a\in L^\vc A^m_{\rho,\delta}$, we have
\begin{equation}\label{eq1}
|\pa_\xi^\alpha a_k(x,y,\xi)|\leq c_\alpha 2^{k(m-\rho|\alpha|)} \ \text{for all $k=1,2,\ldots$}.
\end{equation}

Since $x\in S_j(B), j\geq 2$ and $y, \overline{y}\in B$, we have $x-y\approx x-\overline{y}$. If $|y-\overline{y}|>2^{-k}$, using (a) with $\ell=n+\epsilon$ so that $m-n(\rho-1)-\rho\epsilon+\epsilon<0$, gives
\begin{equation*}
\begin{aligned}
{\rm LHS}&:=\Big|\int a_k(x,y,\xi)e^{i\langle x-y,\xi\rangle}-a_k(x,\overline{y},\xi)e^{i\langle x-\overline{y},\xi\rangle}d\xi\Big|\\
&\leq \Big|\int a_k(x,y,\xi)e^{i\langle x-y,\xi\rangle}d\xi\Big|+\Big|\int a_k(x,\overline{y},\xi)e^{i\langle x-\overline{y},\xi\rangle}d\xi\Big|\\
&\leq C|x-y|^{-n-\epsilon} 2^{k(n+m-\rho n-\rho\epsilon)}\\
&\leq C(2^jr_B)^{-n-\epsilon} 2^{k(m-n(\rho-1)-\rho\epsilon)}.
\end{aligned}
\end{equation*}
This together with the fact that $|y-\overline{y}|>2^{-k}$ gives
\begin{equation}\label{eq1-kernelestimates}
\begin{aligned}
{\rm LHS}&\leq C(2^jr_B)^{-n-\epsilon} 2^{k(m-n(\rho-1)-\rho\epsilon)}\leq C(2^jr_B)^{-n+1} 2^{k((m-n(\rho-1))-\rho\epsilon+\epsilon)}|y-\overline{y}|^\epsilon\\
&\leq C2^{-j\epsilon}(2^jr_B)^{-n} 2^{-k\epsilon'}
\end{aligned}
\end{equation}
where $\epsilon'= -[(m-n(\rho-1))-\rho\epsilon+\epsilon]>0$.

If $|y-\overline{y}|\leq 2^{-k}$, we have
\begin{equation*}
\begin{aligned}
{\rm LHS}&\leq \Big|\int a_k(x,y,\xi)(1-e^{i\langle y-\overline{y},\xi\rangle})e^{i\langle x-y,\xi\rangle}d\xi\Big|\\
&~~~~~+ \Big|\int (a_k(x,y,\xi)-a_k(x,\overline{y},\xi))e^{i\langle x-\overline{y},\xi\rangle}d\xi\Big|: =E_1+E_2.
\end{aligned}
\end{equation*}

We will claim that for all $\ell\geq 0$, we have
\begin{equation}\label{eq1-proff}
E_1\leq C(2^{j}r_B)^\ell 2^{k(m+n-\rho\ell +1)}|y-\overline{y}|.
\end{equation}

Indeed,  we have for all integers $\ell\geq 0$,
\begin{equation*}
\begin{aligned}
E_1&\leq |x-y|^{-\ell }|x-y|^\ell \Big|\int a_k(x,y,\xi)(1-e^{i\langle y-\overline{y},\xi\rangle})e^{i\langle x-y,\xi\rangle}d\xi\Big|\\
&\leq (2^jr_B)^{-\ell }\Big|\sum_{|\alpha|=\ell}\int (x-y)^\alpha a_k(x,y,\xi)\Big(1-e^{i\langle y-\overline{y},\xi\rangle}\Big)e^{i\langle x-y,\xi\rangle}d\xi\Big|\\
&\leq (2^jr_B)^{-\ell }\Big|\sum_{|\alpha|=\ell}\int a_k(x,y,\xi)\Big(1-e^{i\langle y-\overline{y},\xi\rangle}\Big)\pa_\xi^\alpha e^{i\langle x-y,\xi\rangle}d\xi\Big|.
\end{aligned}
\end{equation*}
By integration by part, we get that
\begin{equation}
\begin{aligned}
E_1&\leq (2^jr_B)^{-\ell }\Big|\sum_{|\alpha|=\ell}\int \pa_\xi^\alpha \Big[a_k(x,y,\xi)\Big(1-e^{i\langle y-\overline{y},\xi\rangle}\Big)\Big]e^{i\langle x-y,\xi\rangle}d\xi\Big|.
\end{aligned}
\end{equation}
We write
$$
\sum_{|\alpha|=\ell}\pa_\xi^\alpha \Big[a_k(x,y,\xi)\Big(1-e^{i\langle y-\overline{y},\xi\rangle}\Big)\Big]=\sum_{|\alpha|+|\beta|=\ell} \pa^{\alpha}_\xi a_k(x,y,\xi) \pa^{\beta}_\xi\Big(1-e^{i\langle y-\overline{y},\xi\rangle}\Big).
$$
If $|\beta|=0$, $\Big|1-e^{i\langle y-\overline{y},\xi\rangle}\Big|\leq C|y-\overline{y}||\xi|\leq C2^k|y-\overline{y}|$. Therefore, in this situation,
\begin{equation}
\begin{aligned}
\Big|\sum_{|\alpha|=\ell}\int \pa_\xi^\alpha &\Big[a_k(x,y,\xi)\Big]\Big(1-e^{i\langle y-\overline{y},\xi\rangle}\Big)e^{i\langle x-y,\xi\rangle}d\xi\Big|\\
&\leq C2^{k(n+m+1-\rho|\alpha|)}|y-\overline{y}|=C2^{k(n+m+1-\rho\ell)}|y-\overline{y}|.
\end{aligned}
\end{equation}
Otherwise, $\Big| \pa_\xi^\beta \Big(1-e^{i\langle y-\overline{y},\xi\rangle}\Big)\Big|\leq C|y-\overline{y}|^{|\beta|}$. This together with (\ref{eq1}) gives
\begin{equation}
\begin{aligned}
\Big|\int \pa_\xi^\alpha &a_k(x,y,\xi)\pa_\xi^\beta\Big(1-e^{i\langle y-\overline{y},\xi\rangle}\Big)e^{i\langle x-y,\xi\rangle}d\xi\Big|&\leq C2^{k(n+m-\rho|\alpha|)}|y-\overline{y}|^{|\beta|}\\
&\leq C2^{k(n+m+1-\rho|\alpha|-|\beta|)}|y-\overline{y}|\\
&\leq C2^{k(n+m+1-\rho\ell)}|y-\overline{y}|.
\end{aligned}
\end{equation}
Therefore,
$$
E_1\leq C(2^jr_B)^{-\ell} 2^{k(m+n-\rho\ell+1)}|y-\overline{y}|.
$$
The general statement for non-integer values of $\ell$ follows by interpolation of the inequality for $i$ and $i+1$, where $i<\ell<i+1$. Therefore, (\ref{eq1-proff}) holds for all $\ell>0$. Now taking $\ell = n+ \epsilon$ so that $\epsilon'=-(m+n-\rho n - \rho\epsilon +\epsilon)>0$, we have
$$
\begin{aligned}
E_1&\leq C(2^jr_B)^{-n-\epsilon} 2^{k(m+n-\rho n -\rho\epsilon+\epsilon)}|y-\overline{y}|^\epsilon (2^k|y-\overline{y}|)^{1-\epsilon}\\
&\leq C(2^jr_B)^{-n-\epsilon} 2^{-k\epsilon'}|y-\overline{y}|^\epsilon\\
&\leq C2^{-j\epsilon}(2^jr_B)^{-n} 2^{-k\epsilon'}.
\end{aligned}
$$

It remains to take care the term $E_2$. Repeating the previous arguments we also obtain
\begin{equation*}
\begin{aligned}
E_2&\leq  (2^jr_B)^{-\ell }\Big|\sum_{|\alpha|=\ell}\int \pa_\xi^\alpha \Big[a_k(x,y,\xi)-a_k(x,\overline{y},\xi)\Big] e^{i\langle y-\overline{y},\xi\rangle}d\xi\Big|.
\end{aligned}
\end{equation*}
At this stage, using the Mean value Theorem (apply for each component of $a_k$) and then using the definition of the class $L^\vc A^m_{\rho,\delta}$ give
\begin{equation*}
\begin{aligned}
E_2&\leq  C(2^jr_B)^{-\ell }|y-\overline{y}| 2^{k(n+m-\rho\ell+\delta)}\\
&\leq   C(2^jr_B)^{-\ell }|y-\overline{y}| 2^{k(n+m-\rho\ell+1)}
\end{aligned}
\end{equation*}
for all integer $\ell\geq 0$. Hence, by interpolation again,
$$
E_2\leq   C(2^jr_B)^{-\ell }|y-\overline{y}| 2^{k(n+m-\rho\ell+1)}
$$
for all $\ell\geq 0$. Repeating the arguments used to estimate $E_1$, we conclude that
$$
E_2\leq C2^{-j\epsilon}(2^jr_B)^{-n} 2^{-k\epsilon'}.
$$

Therefore, LHS $\leq C2^{-j\epsilon}(2^jr_B)^{-n} 2^{-k\epsilon'}$. It remains to shows that
\begin{equation}\label{eq2-proof}
{\rm LHS} \leq C2^{-j\epsilon}(2^jr_B)^{-n-N} 2^{-k\epsilon'}.
\end{equation}
To do this, we repeat the arguments above with $\ell = N + n +\epsilon$. Since the proof of this part is analogous to (\ref{eq2-proof}) and hence we omit details here. This completes our proof.

\medskip

(c) If $2^{-k}\leq r_B$, using the argument as in (b), we have
\begin{equation*}
\begin{aligned}
{\rm LHS}&:=\Big|\int a_k(x,y,\xi)e^{i\langle x-y,\xi\rangle}-a_k(x,\overline{y},\xi)e^{i\langle x-\overline{y},\xi\rangle}d\xi\Big|\\
&\leq \Big|\int a_k(x,y,\xi)e^{i\langle x-y,\xi\rangle}d\xi\Big|+\Big|\int a_k(x,\overline{y},\xi)e^{i\langle x-\overline{y},\xi\rangle}d\xi\Big|\\
&\leq C|x-y|^{-n-\epsilon} 2^{-k\epsilon}\\
&\leq C\f{(2^jr_B)^{-n-\epsilon}}{r_B^\epsilon}\f{2^{-k\epsilon}}{r_B^\epsilon}=C2^{-j\epsilon}(2^jr_B)^n\Big(\f{1}{r_B2^k}\Big)^\epsilon.
\end{aligned}
\end{equation*}

If $r_B< 2^{-k}$, we have
\begin{equation*}
\begin{aligned}
{\rm LHS}&\leq \Big|\int a_k(x,y,\xi)(1-e^{i\langle y-\overline{y},\xi\rangle})e^{i\langle x-y,\xi\rangle}d\xi\Big|\\
&~~~~~+ \Big|\int (a_k(x,y,\xi)-a_k(x,\overline{y},\xi))e^{i\langle x-\overline{y},\xi\rangle}d\xi\Big|: =E_1+E_2.
\end{aligned}
\end{equation*}
The previous arguments in (b) show that
$$
\begin{aligned}
E_1+E_2&\leq C(2^jr_B)^{-n-\epsilon} 2^{k(-\epsilon +1)}|y-\overline{y}|\\
&\leq C(2^jr_B)^{-n-\epsilon} 2^{k(-\epsilon +1)}r_B=C(2^jr_B)^{-n-\epsilon}r_B^\epsilon (r_B2^{k})^{(-\epsilon +1)}\\
&\leq C2^{-j\epsilon}(2^jr_B)^{-n}(r_B2^{k})^{(1-\epsilon)}.
\end{aligned}
$$
Hence,
$$
\Big|\int a_k(x,y,\xi)e^{i\langle x-y,\xi\rangle}-a_k(x,\overline{y},\xi)e^{i\langle x-\overline{y},\xi\rangle}d\xi\Big|\leq C2^{-j\epsilon}(2^jr_B)^{-n}(2^{k}r_B)^{\epsilon'} \ \text{if $2^{k}r_B\leq 1$}
$$
and
$$
\Big|\int a_k(x,y,\xi)e^{i\langle x-y,\xi\rangle}-a_k(x,\overline{y},\xi)e^{i\langle x-\overline{y},\xi\rangle}d\xi\Big|\leq C2^{-j\epsilon}(2^jr_B)^{-n}(2^{k}r_B)^{-\epsilon'} \ \text{if $2^{k}r_B> 1$}.
$$
By taking $\ell = n+N+\epsilon$ and repeating the previous arguments, we obtain
$$
\Big|\int a_k(x,y,\xi)e^{i\langle x-y,\xi\rangle}-a_k(x,\overline{y},\xi)e^{i\langle x-\overline{y},\xi\rangle}d\xi\Big|\leq C2^{-j\epsilon}(2^jr_B)^{-n-N}(2^{k}r_B)^{\epsilon'} \ \text{if $2^{k}r_B\leq 1$}
$$
and
$$
\Big|\int a_k(x,y,\xi)e^{i\langle x-y,\xi\rangle}-a_k(x,\overline{y},\xi)e^{i\langle x-\overline{y},\xi\rangle}d\xi\Big|\leq C2^{-j\epsilon}(2^jr_B)^{-n-N}(2^{k}r_B)^{-\epsilon'} \ \text{if $2^{k}r_B> 1$}.
$$
This completes the proof of (c).

\medskip

Since the associated kernel $K(x,y)$ of the operator $T_a$ is given by $$
K(x,y)=\f{1}{(2\pi)^n}\int a(x,y,\xi)e^{i\langle x-y, \xi\rangle}d\xi = \sum_{k\geq 0}\f{1}{(2\pi)^n}\int a_k(x,y,\xi)e^{i\langle x-y, \xi\rangle}d\xi
$$
with $a_k(x,\xi)$ as in Lemma \ref{lem1}, from Lemma \ref{lem1} we imply directly the following result.

\begin{lem}\label{lem2}
Let $a\in L^\vc A^m_{\rho,\delta}$ with $m<n(\rho-1)$ or $a\in L^\vc A^{0}_{1,\delta}, \delta\in [0,1]$ and let $K^*(x,y)$ be the associated kernel of the operator $T_a^*$, the conjugate of $T_a$.

(a) For any $N >0$, we have
$$
|K^*(x,y)|\leq \f{C}{|x-y|^{-N}}, \ x\neq y;
$$
\medskip
(b) For any $N>0$, there exists $\epsilon >0$ so that any ball $B\subset \RR^n$, $y, \overline{y} \in B, x \in S_j(B), j\geq 2$, we have
$$
|K^*(y,x)-K^*(\overline{y},x)|\leq C2^{-j\epsilon}(2^jr_B)^{-n}\min\{1, (2^jr_B)^{-N}\}.
$$
\end{lem}

\section{Proof of Theorem \ref{mainresult}}
Note that the boundedness of $T_a$ and $[b, T_a]$ on $L^p(w)$, $1<p<\vc$ and $w\in A^\vc_p$ is equivalent to that of $T_a^*$. Therefore, it is suffices to prove (a) and (b) for $T^*_a$.
For $b\in BMO_\vc$, set $T^{*,b}_a=[b,T^*_a]$. Before coming to the proof of Theorem \ref{mainresult}, we need the following results.
\begin{lem}\label{lem1-thm1}
Let $a\in L^\vc A^m_{\rho,\delta}$ with $m<n(\rho-1)$ or $a\in L^\vc A^0_{1,\delta}, \delta\in [0,1]$ and $b\in BMO_\theta, \theta\geq 0$. If $T_a$ is bounded on $L^p$ for all $1<p<\vc$, then for any $p>1$ there exists $C>0$ such that for all balls $Q=Q(x_0,1)$,

(a)
$$
\f{1}{|Q|}\int_{Q}|T^*_af(x)|dx\leq C\inf_{y\in Q}G_{p}(y);
$$

(b)
$$
\f{1}{|Q|}\int_{Q}|T^{*,b}_af(x)dx|\leq C\inf_{y\in Q}G_{p}(y)\|b\|_\theta;
$$
\end{lem}
\emph{Proof:}

(a) We split $f=f_1+f_2$ where $f_1=f\chi_{4Q}$. For each $j\geq 0$, we have
$$
\f{1}{|Q|}\int_{Q}|T^*_af(x)|dx \leq \f{1}{|Q|}\int_{Q}|T^*_af_1| +\f{1}{|Q|}\int_{Q}|T^*_af_2|:=I_{1}+I_{2}.
$$
Using H\"older inequality and the fact that $T_a^*$ is bounded on $L^p, 1<p<\vc$, we write
\begin{equation*}
\begin{aligned}
I_{1}&\leq  C \Big(\f{1}{|Q|}\int_{Q}|T_a^*f_1|^{p}\Big)^{1/p}\leq \Big(\f{1}{|4Q|}\int_{4Q}|f|^{p}\Big)^{1/p}\\
& \leq  C \inf_{y\in Q}G_{p}f(y).
\end{aligned}
\end{equation*}
For the term $I_{2}$ we have, for $x\in Q$,
\begin{equation*}
\begin{aligned}
T_a^*f_2(x)&=\int_{R^n\backslash 4Q}K^*(x,y)f(y)dy=\int_{R^n\backslash 4Q}K^*(y,x)f(y)dy\\
& =\sum_{k\geq 3}\int_{S_k(Q)}K^*(y,x)f(y)dy.
\end{aligned}
\end{equation*}
Applying (a) of Lemma \ref{lem2}, we have
\begin{equation}\label{eq1-prooflem1-mainthm}
\begin{aligned}
T_a^*f_2(x)=\sum_{k\geq 3}\int_{S_k(Q)}K^*(x,y)f(y)dy&\leq \sum_{k\geq 3}\int_{S_k(Q)}\f{f(y)}{|x-y|^{n+N}}dy\\
&\leq C\inf_{y\in Q}Gf(y)\leq C\inf_{y\in Q}G_pf(y).
\end{aligned}
\end{equation}
This completes the proof of (a).

\medskip

(b)  Taking $1<r<p$, we write
$$
T_a^{*,b} f=(b-b_Q)T_a^* f -T_a^*((b-b_Q)f).
$$
So, we have
\begin{equation*}
\begin{aligned}
\f{1}{|Q|}\int_{Q}|T_a^{*,b} f(x)|d&\leq \f{1}{|Q|}\int_{Q}|(b-b_Q)T_a^* f|dx+\f{1}{|Q|}\int_{Q}|T_a^*((b-b_Q)f)(x)|dx\\
&: =I_1+I_2.
\end{aligned}
\end{equation*}
Let us estimate $I_1$ first. By H\"older inequality, we can write
\begin{equation*}
\begin{aligned}
I_{1}&\leq  C \|b\|_{\theta} \Big(\f{1}{|Q|}\int_{Q}|T_a^* f|^{p}\Big)^{1/p}\\
& \leq  C \|b\|_{\theta}\Big(\Big(\f{1}{|Q|}\int_{Q}|T_a^* f_1|^{p}\Big)^{1/p}+\Big(\f{1}{|Q|}\int_{Q}|T_a^* f_2|^{p}\Big)^{1/p}\Big)\\
&:=I_{11}+I_{12}
\end{aligned}
\end{equation*}
where $f=f_1+f_2$ with $f_1=f\chi_{4Q}$.

Due to $L^{p}$-boundedness of $T_a^*$, one has
\begin{equation*}
\begin{aligned}
I_{11}\leq  C \Big(\f{1}{|4Q|}\int_{4Q}|f|^{p}\Big)^{1/p}\leq  C \inf_{y\in Q}G_{p}f(y).
\end{aligned}
\end{equation*}

To estimate $I_{12}$, using (\ref{eq1-prooflem1-mainthm}) gives
$I_{12}\leq C\inf_{y\in Q}G_{p}f(y)$.

\medskip

The estimate for $I_2$ can be proceeded in the same method. Indeed, we write
\begin{equation*}
\begin{aligned}
\f{1}{|Q|}\int_{Q}&|T_a^*((b-b_Q)f)(x)|dx\\
&\leq \f{1}{|Q|}\int_{Q}|T_a^*((b-b_Q)f_1)(x)|dx+\f{1}{|Q|}\int_{Q}|T_a^*((b-b_Q)f_2)(x)|dx\\
&:=I_{21}+I_{22}
\end{aligned}
\end{equation*}
where $f=f_1+f_2$ and $f_1=f\chi_{4Q}$.

To estimate $I_{21}$, using H\"older inequality we have
\begin{equation*}
\begin{aligned}
\f{1}{|Q|}\int_{Q}&|T_a^*((b-b_Q)f_1)(x)|dx\\
 &\leq \Big(\f{1}{|Q|}\int_{Q}|T_a^*((b-b_Q)f_1)(x)|^rdx\Big)^{1/r}\\
&\leq \Big(\f{1}{|Q|}\int_{Q}|((b-b_Q)f_1)(x)|^rdx\Big)^{1/r}\\
&\leq \Big(\f{1}{|4Q|}\int_{4Q}|f(x)|^{p}dx\Big)^{1/p}\Big(\f{1}{|4Q|}\int_{4Q}|b(x)-b_Q|^{\nu}dx\Big)^{1/\nu} \ \ (\nu=\f{pr}{p-r})\\
&\leq C\|b\|_{\theta}\inf_{y\in Q}G_{p}f(y).
\end{aligned}
\end{equation*}

For the term $I_{22}$, due to (a) of Lemma \ref{lem2}, we can write
\begin{equation}\label{eq1-lem1-thm2}
\begin{aligned}
T_a^* ((b-b_Q)f_2)(x)
&=\sum_{k\geq 3}\int_{S_k(Q)}K^*(x,y)((b-b_Q)f)(y)dy\\
&\leq C\sum_{k\geq 3}2^{-kN}\int_{S_k(Q)}|(b(y)-b_Q)f(y)|dy.
\end{aligned}
\end{equation}

By H\"older inequality and Proposition \ref{JNforBMOL}, we have give
\begin{equation}\label{eq3-lem1-thm2}
\begin{aligned}
\int_{S_k(Q)}&|(b(y)-b_Q)f(y)|dy\\
&\leq |2^kQ|\Big(\f{1}{|2^kQ|}\int_{2^kQ}|f|^{p}\Big)^{1/p}\Big(\f{1}{|2^kQ|}\int_{2^kQ}|b-b_Q|^{p'}\Big)^{1/p'} \\
&\leq k2^{k\theta}|2^kQ|\|b\|_{\theta}\Big(\f{1}{|2^kQ|}\int_{2^kQ}|f|^{p}\Big)^{1/p}\\
&\leq k2^{k(\theta+n)}\|b\|_{\theta}\Big(\f{1}{|2^kQ|}\int_{2^kQ}|f|^{p}\Big)^{1/p}.
\end{aligned}
\end{equation}

From  (\ref{eq3-lem1-thm2}) and (\ref{eq1-lem1-thm2}) we obtain
\begin{equation}\label{eq1-lem1-thm2}
\begin{aligned}
T_a^* &((b-b_Q)f_2)(x) \leq C\|b\|_{\theta}\inf_{y\in Q}G_pf(y).
\end{aligned}
\end{equation}

This completes our proof.

\begin{rem}\label{rem2}
{\rm The result in Lemma \ref{lem1-thm1} still holds if we replace the critical ball $Q$ by $2Q$}.
\end{rem}

\begin{lem}\label{lem2-thm1}
Let $a\in L^\vc A^m_{\rho,\delta}$ with $m<n(\rho-1)$ or $a\in L^\vc A^0_{1,\delta}, \delta\in [0,1]$ and $b\in BMO_\theta, \theta\geq 0$. If $T_a$ is bounded on $L^p$ for all $1<p<\vc$, then for any $p>1$ there exists $C>0$ so that for all $f$ and $x,y\in B=B(x_B,r_B)$ with $r_B<4$, we have

(a)
$$
\int_{\mathbb{R}^n\backslash 2B}|(K^*(x,z)-K^*(y,z))f(z)|dz\leq C(\inf_{u\in B}G_{4,p}f(u)+\inf_{u\in B}\widetilde{M}_{p}f(u));
$$

(b)
$$
\int_{\mathbb{R}^n\backslash 2B}|(K^*(x,z)-K^*(y,z))((b-b_B)f)(z)|dz\leq C\|b\|_{\theta}\inf_{u\in B}(\inf_{u\in B}G_{4,p}f(u)+\inf_{u\in B}\widetilde{M}_{p}f(u)).
$$

\end{lem}
\emph{Proof:}

(a) Using (b) of Lemma \ref{lem2}, we write
\begin{equation*}
\begin{aligned}
\int_{\mathbb{R}^n\backslash 2B}&|(K^*(x,z)-K^*(y,z))f(z)|dz\\
&\leq C\sum_{k\geq 2} \int_{S_k(B)}|(K^*(x,z)-K^*(y,z))f(z)|dz\\
&\leq C\sum_{k\geq 2}2^{-k\epsilon}(2^kr_B)^{-n}\min\{1, (2^jr_B)^{-N}\}\int_{S_k(B)}|f(z)|dz\\
&\leq C\sum_{k\geq 2}2^{-k\epsilon}\min\{1, (2^jr_B)^{-N}\}\f{1}{|2^kB|}\int_{S_k(B)}|f(z)|dz\\
&= \sum_{k= 2}^{k_0}\ldots+\sum_{k> k_0}\ldots:= I_1+I_2
\end{aligned}
\end{equation*}
where $k_0$ is the smallest integer so that $2^{k_0+1}r_B>4$.

To estimate $I_1$, Let $\{Q_l\}$ and $\{\widetilde{Q}_l\}$ be families of balls as in (\ref{defnoftidleM}). If $x\in Q_l\cap B$ then $2^kB\subset \widetilde{Q}_l$ for all $k=1, 2, \ldots k_0$. This implies
$$
\f{1}{|2^kB|}\int_{2^kB}|f(z)|dz\leq \inf_{u\in B}\widetilde{M}_{p}f(u)
$$
for all $k=1, 2, \ldots k_0$.

Hence
\begin{equation*}
\begin{aligned}
I_1\leq \sum_{k= 2}^{k_0} 2^{-k\epsilon}\inf_{u\in B}\widetilde{M}_{p}f(u)\leq C\inf_{u\in B}\widetilde{M}_{p}f(u).
\end{aligned}
\end{equation*}

For the term $I_2$, since $2^{k_0}r_B\geq 4$ we have
\begin{equation*}
\begin{aligned}
I_2&\leq \sum_{k\geq k_0} 2^{-k\epsilon}(2^kr_B)^{-N}\f{1}{|2^kB|}\int_{S_k(B)}|f(z)|dz\\
&\leq \sum_{k\geq k_0} 2^{-k\epsilon}(2^{k-k_0} 2^{k_0}r_B)^{-N}\f{1}{|2^{k-k_0} 2^{k_0}B|}\int_{2^{k-k_0} 2^{k_0}B}|f(z)|dz\\
&\leq \sum_{k\geq k_0} 2^{-k\epsilon}(2^{k-k_0})^{-N}\f{1}{|2^{k-k_0} 2^{k_0}B|}\int_{2^{k-k_0} 2^{k_0}B}|f(z)|dz\\
&\leq \sum_{k\geq 0} 2^{-k\epsilon}2^{-kN}\f{1}{|2^{k} 2^{k_0}B|}\int_{2^{k} 2^{k_0}B}|f(z)|dz.
\end{aligned}
\end{equation*}
Note that $2^{k_0}B\subset \widehat{Q}=4 Q$ here $Q=B(x_0,1)$ and $|Q|\approx |2^{k_0}B|$. So, we have
\begin{equation*}
\begin{aligned}
I_2&\leq \sum_{k\geq 0} 2^{-k\epsilon}2^{-kN} \Big(\f{1}{|2^{k}\widehat{Q}|}
\int_{2^{k}\widehat{Q}}|f(z)|dz\Big)\\
&\leq C\inf_{u\in B}G_{4,p}f(u).
\end{aligned}
\end{equation*}
Hence, we get (a).

\medskip

(b) Using H\"older inequality and (b) of Lemma \ref{lem2}, we obtain
\begin{equation*}
\begin{aligned}
\int_{\mathbb{R}^n\backslash 2B}&|(K^*(x,z)-K^*(y,z))((b-b_B)f)(z)|dz\\
&= \sum_{k\geq 2} \int_{S_k(B)}|(K^*(x,z)-K^*(y,z))((b-b_B)f)(z)|dz\\
&\leq C\sum_{k\geq 2}2^{-k\epsilon}\min\{1, (2^jr_B)^{-N}\}\f{1}{|2^kB|}\int_{S_k(B)}|((b-b_B)f)(z)|dz\\
&\leq C\sum_{k\geq 2}2^{-k\epsilon}\min\{1, (2^jr_B)^{-N}\} \Big(\f{1}{|2^kB|}\int_{2^kB}|f(z)|^{p}dz\Big)^{1/p}\Big(\f{1}{|2^kB|}\int_{2^kB}|b(z)-b_B|^{p'}dz\Big)^{1/p'}.
\end{aligned}
\end{equation*}

Now using Proposition \ref{JNforBMOL}, we get that
\begin{equation*}
\begin{aligned}
\int_{\mathbb{R}^n\backslash 2B}&|(K^*(x,z)-K^*(y,z))((b-b_B)f)(z)|dz\\
&\leq C\sum_{k\geq 2}k2^{-k\epsilon}(2^jr_B)^\theta \min\{1, (2^jr_B)^{-N}\}\|b\|_{\theta} \Big(\f{1}{|2^kB|}\int_{2^kB}|f(z)|^{p}dz\Big)^{1/p}.
\end{aligned}
\end{equation*}

At this stage, repeating the same argument as in (a), we complete the proof of (b).

\medskip

We are now in position to prove Theorem \ref{mainresult}.

\medskip

\emph{Proof of Theorem \ref{mainresult}:} (a) Using the standard argument, see for example \cite{BHS1}, fix $1<p<\vc$ and $w\in A_p^\vc$. So, by Proposition \ref{property-Avc} we can pick $r>1$ and $\nu\geq 0$ so that $w\in A^\nu_{p/r}$. By Proposition \ref{FSinequalityversion} we have
\begin{equation*}
\begin{aligned}
\|T_a^*f\|^p_{L^p(w)}&\leq \|M_{{\rm loc},\beta}T_a^*f\|^p_{L^p(w)}\\
&\leq C\|M^\sharp_{{\rm loc},4}T_a^*f\|^p_{L^p(w)} + C\sum_{k}w(Q_k)\Big(\f{1}{2Q_k}\int_{2Q_k}|T_a^*f|\Big)^p\\
&:=I_1+I_2.
\end{aligned}
\end{equation*}
Let us estimate $I_1$ first. By Lemma \ref{lem1-thm1} and Remark \ref{rem2}, we have
$$
\f{1}{2Q_k}\int_{2Q_k}|T_a^*f|\leq C\inf_{y\in Q_k}G_{r}f(y).
$$

Invoking Proposition \ref{weighted estiamtes for G}, we conclude that
\begin{equation}\label{eq1-proofthm1}
\begin{aligned}
\sum_{k}w(Q_k)\Big(\f{1}{2Q_k}\int_{2Q_k}|T_a^*f|\Big)^p&\leq \sum_k \int_{Q_k}|G_{r}f(x)|^{p}w(x)dx\\
&\leq C\int_{\mathbb{R}^n}|G_{r}f(x)|^{p}w(x)dx\\
&\leq C\|f\|^p_{L^p(w)}.
\end{aligned}
\end{equation}

We now take care $I_2$. For any ball $B(x_0, r_B)$ with $r_B\leq 4$ and $x\in B$, we write
\begin{equation*}
\begin{aligned}
\f{1}{|B|}\int_B&|T_a^*f(x)-(T_a^*f)_B|dx\\
&\leq \f{2}{|B|}\int_B|T_a^*f_1(x)|dx +\f{1}{|B|}\int_B|T_a^*f_2(x)-(T_a^*)f_2))_B|dx\\
&:=E_1+E_2.
\end{aligned}
\end{equation*}
where $f=f_1+f_2$ with $f_1=f\chi_{2B}$.

For $E_1$, since $T_a^*$ is bounded on $L^{r}$, we have
\begin{equation*}
\begin{aligned}
\f{1}{|B|}\int_B|T_a^*f_1(x)|dx&\leq \Big(\f{1}{|B|}\int_B|T_a^*f_1(x)|^{r}dx\Big)^{r}\\
&\leq C\Big(\f{1}{|2B|}\int_{2B}|f|^{r}dx\Big)^{r}\\
&\leq C\inf_{u\in B}\widetilde{M}_{r}f(u).
\end{aligned}
\end{equation*}
Due to Lemma \ref{lem2-thm1}, we write
\begin{equation*}
\begin{aligned}
E_2&\leq \f{1}{|B|^2}\int_B\int_B\Big(\int_{\mathbb{R}^n\backslash 2B}|(K^*(u,z)-K^*(y,z))f(z)|dz\Big)dydu\\
&\leq C(\inf_{u\in B}G_{4,r}f(u)+\inf_{u\in B}\widetilde{M}_{r}f(u)).
\end{aligned}
\end{equation*}
These two estimates of $E_1$ and $E_2$ tell us that
$$
M^\sharp_{{\rm loc},4}T_a^*f(x)\leq C(G_{4,r}f(x)+\widetilde{M}_{r}f(x)).
$$
Applying Proposition \ref{weighted estiamtes for G} and the weighted estimates of $\widetilde{M}_{r}$, we get that
\begin{equation}\label{eq2-proofthm1}
\|M^\sharp_{{\rm loc},4}T_a^*f\|_{L^p(w)}\leq C\|f\|_{L^p(w)}.
\end{equation}
From (\ref{eq1-proofthm1}) and (\ref{eq2-proofthm1}), we obtain
$$
\|T_a^*f\|_{L^p(w)}\leq C\|f\|_{L^p(w)}.
$$
This completes our proof.

\medskip

(b) Fixed $1<p<\vc, b\in BMO_\theta, \theta\geq 0$ and $w\in A_p^\vc$. So, we can pick $r>1$ and $\nu\geq 0$ so that $w\in A^\nu_{p/r}$. Then we have by Lemma \ref{FSinequalityversion}
\begin{equation*}
\begin{aligned}
\|T_a^{*,b}f\|^p_{L^p(w)}&\leq \int_{\mathbb{R}^n}|M_{{\rm loc},\beta}(T_a^{*,b}f)(x)|^pw(x)dx\\
&\leq C \int_{\mathbb{R}^n}|M^\sharp_{{\rm loc},4}(T_a^{*,b}f)(x)|^pw(x)dx+\sum_{k}w(Q_k)\Big(\f{1}{|2Q_k|}\int_{2Q_k}|T_a^{*,b}f|\Big)^p
\end{aligned}
\end{equation*}
where $\{Q_k\}$ is a family of critical balls given in Lemma \ref{FSinequalityversion}.

The analogous argument to that in (a) gives
\begin{equation*}
\begin{aligned}
\sum_{k}w(Q_k)\Big(\f{1}{|2Q_k|}&\int_{2Q_k}|T_a^{*,b}f|\Big)^p\leq C\|b\|^p_\theta\|f\|_{L^p(w)}^p.
\end{aligned}
\end{equation*}

It remains to estimate $\int_{\mathbb{R}^n}|M^\sharp_{{\rm loc},4}(T_a^{*,b}f)(x)|^pw(x)dx$. For any ball $B(x_0, r_B)$ with $r_B\leq 4$ and $x\in B$, we write
\begin{equation*}
\begin{aligned}
\f{1}{|B|}\int_B&|T_a^{*,b}f(x)-(T_a^{*,b}f)_B|dx\\
&\leq \f{2}{|B|}\int_B|(b-b_B)T^*_af(x)|dx+\f{2}{|B|}\int_B|T^*_a((b-b_B)f_1)(x)|dx\\
& ~~~~~ +\f{1}{|B|}\int_B|T^*_a((b-b_B)f_2)(x)-(T^*_a((b-b_B)f_2))_B|dx\\
&:=E_1+E_2+E_3.
\end{aligned}
\end{equation*}
where $f=f_1+f_2$ with $f_1=f\chi_{2B}$.

H\"older inequality and Proposition \ref{JNforBMOL} show that
\begin{equation*}
\begin{aligned}
 E_1&\leq C\Big(\f{1}{|B|}\int_B|b-b_B|^{r'}\Big)^{1/r'}\Big(\f{1}{|B|}\int_B|T^*_a|^{r}\Big)^{1/r}\\
 &\leq C\|b\|_{\theta} \Big(\f{1}{|B|}\int_B|T^*_af|^{r}\Big)^{1/r}.
\end{aligned}
\end{equation*}

For any critical ball $Q_j$ such that $x\in Q_j\cap B$. It can be verified that $B\subset \widetilde{Q}_j:=8Q_j$. This yields that
$$
E_1 \leq C\|b\|_{\theta}\times \inf_{y\in B}\widetilde{M}_{r}(T^*_af)(y).
$$
Using H\"older inequality and and Proposition \ref{JNforBMOL} again, we have, for $1<s<r$,
\begin{equation*}
\begin{aligned}
 E_2&\leq C\Big(\f{1}{|B|}\int_B|T^*_a((b-b_B)f_1)|^{s}\Big)^{1/s}\\
 &\leq C\Big(\f{1}{|B|}\int_{2B}|(b-b_B)f_1|^{s}\Big)^{1/s}\\
&\lesssim \Big(\f{1}{|B|}\int_{2B}|(b-b_B)|^{\gamma}\Big)^{1/\gamma}\Big(\f{1}{|B|}\int_{2B}|f|^{r}\Big)^{1/r} \ \text{for some $\gamma>s$}\\
&\lesssim \|b\|_{\theta}\times \inf_{y\in B}\widetilde{M}_{r}(f)(y).
\end{aligned}
\end{equation*}

To estimate $E_3$, using Lemma \ref{lem2-thm1}, we conclude
\begin{equation*}
\begin{aligned}
E_3&\leq C\f{1}{|B|^2}\int_B\int_B\Big(\int_{R^n\backslash 2B}|K^*(u,z)-K^*(y,z)||b(z)-b_B||f(z)|dz\Big)dydu\\
&\leq C\|b\|_\theta(G_{4, r}f(x)+\widetilde{M}_{r}(f)(x)).
\end{aligned}
\end{equation*}
These three estimates of $E_1, E_2$ and $E_3$ give
$$
M^\sharp_{{\rm loc},4}(T^{*,b}_af)(x)\leq C\|b\|_\theta(\widetilde{M}_{r}(T^*_af)(x)+G_{4, r}(x)+\widetilde{M}_{r}(f)(x)).
$$
This implies
$$
\|M^\sharp_{{\rm loc},4}(T^{*,b}_af)\|_{L^p(w)}\leq C\|b\|_\theta(\|\widetilde{M}_{p_0}(T^*_af)\|_{L^p(w)}+\|G_{4, r}f\|_{L^p(w)}+\|\widetilde{M}_{r}(f)\|_{L^p(w)}.
$$
Since $\widetilde{M}_{r}, G_{4, r}$ and $T^*_a$ is bounded on $L^p(w)$, we obtain the desired results.\\

This completes our proof.\\

\end{document}